\pdfoutput=1

\documentclass[final,5p,times,twocolumn]{elsarticle}

\makeatletter
\def\ps@pprintTitle{%
 \let\@oddhead\@empty
 \let\@evenhead\@empty
 \def\@oddfoot{}%
 \let\@evenfoot\@oddfoot}
\makeatother

\usepackage{amsmath,amssymb}

\usepackage[justification=centering]{caption}
\usepackage{subcaption}
\usepackage{float}
\usepackage{color}

\newtheorem{proposition}{Proposition}
\newtheorem{problem}{Problem}

\journal{}

\begin{document}

\begin{frontmatter}

\title{Approximate explicit model predictive control via piecewise nonlinear system identification\tnoteref{label0}}
\tnotetext[label0]{This work is partially supported by the French National Research Agency, task agreement ANR-13-SEED-0005.}

\author[cnrs,cea]{Van-Vuong~Trinh\corref{cor1}}
\cortext[cor1]{Corresponding author}
\ead{van-vuong.trinh@gipsa-lab.fr}
\ead[url]{http://trinhvv.github.io}

\author[cnrs]{Mazen~Alamir}
\ead{mazen.alamir@gipsa-lab.fr} 

\author[cea]{Patrick~Bonnay}
\ead{patrick.bonnay@cea.fr}

\address[cnrs]{Univ. Grenoble Alpes, CNRS, GIPSA-lab, F-38000 Grenoble, France
} 
\address[cea]{Univ. Grenoble Alpes, CEA, INAC-SBT, F-38000 Grenoble, France}    


\begin{abstract}
This article presents an identification methodology to capture general relationships, with application to piecewise nonlinear approximations of model predictive control for constrained (non)linear systems. The mathematical formulation takes, at each iteration, the form of a constrained linear (or quadratic) optimization problem that is mathematically feasible as well as numerically tractable. The efficiency of the devised methodology is demonstrated via two industrial applications. Results suggest the possibility to achieve high approximate precision with limited number of regions, leading to a significant reduction in computation time when compared to the state-of-the-art implicit model predictive control solvers.
\end{abstract}

\end{frontmatter}

\section{Introduction}

Model predictive control (MPC) \cite{rawlings2009model} is a model-based control method which involves the solution at each sampling period of a finite horizon optimal control problem subject to the system dynamics and constraints. This form of MPC, whose control actions are obtained by an online optimization, is referred to as implicit MPC. To reduce the online computational complexity of such an optimal control problem, two directions have been investigated. On the one hand, more efficient optimization strategies have been proposed in the aforementioned implicit MPC context. These approaches try to exploit as far as possible the problem structure and the possibility to use hot starts in order to reduce the number of iterations needed to reach a good solution (see e.g. \cite{diehl2009efficient} and the references therein). On the other hand, the so-called explicit MPC (EMPC) has been the subject of significant research efforts, where the control laws are available a priori, accordingly reducing an online optimization to a function evaluation. This contribution lies in this second direction with nonlinear systems as targeted domain of application.

For linear or hybrid systems, the explicit control laws can be exactly solved offline via parametric programming \cite{bemporad2002explicit} \cite{oberdieck2015explicit}. As for nonlinear systems, deriving the true optimal nonlinear MPC (NMPC) control law is generally not possible and hence approximate approaches have to be followed. Following this direction, several techniques have been devised where the MPC control laws are approximated using piecewise affine approximators \cite{grancharova2012explicit}, artificial neural networks \cite{aakesson2006neural} \cite{pin2013approximate} and set memberships identification \cite{canale2010efficient} \cite{fagiano2012set}. Nevertheless, piecewise affine approaches normally lead to an excessive number of regions needed to represent the approximate control law while standard nonlinear approximators offer universal capabilities at the price of non-convex optimization schemes. Thus, approximating (or learning) the MPC control laws by piecewise nonlinear (PWNL) functions are recently of interest (see, for instance, \cite{domahidi2011learning}).

This paper suggests a practical computational methodology to EMPCs derived from learning data that incorporates implicit MPC solution and a problem-dependent regressor. In particular, each component of the control vector is identified independently as a PWNL map of the available regressor. A heuristic identification procedure to do so is proposed, which is an extension of the preliminary work \cite{alamir2013new}. Regarding complexity, the PWNL representations of the identified EMPCs allows to reduce the number of regions with respect to the piecewise affine approaches, especially for constrained NMPCs. Compared to the standard nonlinear structures, the methodology takes advantages of efficient computation of constrained linear (or quadratic) programming problems. 

The preliminary version of this paper has been presented in \cite{trinh2016explicit}. The current work introduces additional details, including the improvement of algorithms, their numerical analysis as well as one more multivariable control application, namely, the constrained control of the compression station used in the cryogenic refrigerators.

The remainder of this paper is organized as follows. Section \ref{sec:nl} recalls the identification of a nonlinear map proposed in \cite{alamir2013new}. The extension of this framework to the piecewise nonlinear maps is explained in Section \ref{sec:pwnl} where an identification methodology is developed with application to EMPCs. Two real control design problems are used in Section \ref{sec:apps} to assess the ability of the proposed methodology to approximate the MPC control laws with rather small number of regions and the possibility to monitor the complexity as a function of the required precision. Moreover, the resulting computation time is compared to the one obtained by the state-of-the-art solver ACADO \cite{houska2011acado}, showing the potential benefit from the explicit representation to extend the applicability of the MPC framework. Finally, Section \ref{sec:conclusion} concludes the paper and gives hints for further investigations.

\subsection*{Notation}

In the sequel, the following notation is extensively used: For $x \in \mathbb{R}^n$, its $i$-th component is $x_i$, the Euclidean norm is $||x||$, the weighted norm is $||x||_A=(x^TAx)^{1 {}/ 2}$, the $p$-norm is $||x||_p$ with $p \geq 1$ be a real number, the absolute norm is $|x|$ when $n=1$. For a discrete set, $\#\mathcal{X}$ means the cardinality, i.e. the number of elements of the set $\mathcal{X}$.

\section{Design of A Nonlinear Approximator}\label{sec:nl}

In this section, we briefly recall the identification framework \cite{alamir2013new} for a class of nonlinear relationships. Several modifications with preliminary analysis are also introduced.

\subsection{Problem Statement and Mathematical Formulation}

\begin{problem}\label{pb:single_id}
Given the data $\mathcal{D}=\{(q(k),Z(k)\}_{k=1}^N$ where $Z \in \mathcal{Z} \subseteq \mathbb{R}^{n_z}$ is the regressor and $q \in [\underline{q}, \overline{q}] \subset \mathbb{R}$ is the output. Find a map $F : \mathbb{R}^{n_z}  \rightarrow \mathbb{R}$ of the form: 
\begin{equation}\label{eq:approx_struct}
 	F(Z):=\Gamma^{-1}(L^TZ)
\end{equation}
where $L \in \mathbb{R}^{n_z}$, $\Gamma(\cdot)$ is strictly increasing and such that the following approximation holds:
\begin{equation}\label{eq:approx_map}
 	q \approx F(Z)
\end{equation}    
\end{problem}

The structure \eqref{eq:approx_struct} is a Wiener-like model with a strictly increasing static output mapping whose inverse map can be parametrized  using a finite function basis:
\begin{equation}\label{eq:fcn_param}
 	\Gamma(q)=\sum_{i=1}^{n_b}B^{(i)}(\xi(q))\mu_i=B(\xi(q))\mu \text{; } \xi(q)=\frac{q-\underline{q}}{\overline{q}-\underline{q}}
\end{equation} 
where the basis functions are given by: 
\begin{equation}\label{eq:basis_fcn}
\{B^{(i)}\}_{i = 1}^{n_b}:=\{1\} \cup \{B_1^{(i)}\}_{i = 2}^{n_m}\cup \{B_2^{(i)}\}_{i = 1}^{n_m}
\end{equation}
\begin{equation}
B_1^{(i)}(\eta):=(1+\alpha_i) \frac{\eta}{1+\alpha_i \eta} \text{; } B_2^{(i)}(\eta):=\frac{\eta}{1+\alpha_i (1-\eta)}
\end{equation}
The number of functions is $n_b=2n_m$ while the coefficients $\alpha_i$ are given by $\alpha_i:=e^{\beta (1-i)}-1$ for some constant $\beta > 0$. 

Denote $\mu \in \mathbb{R}^{n_b}$ and $L \in \mathbb{R}^{n_z}$ as the parameters of $F(Z)$. The total number of parameters is $n_p=n_b+n_z$. The basic idea of the above formulation is to solve  Problem~\ref{pb:single_id} by finding $\mu \in \mathbb{R}^{n_b}$ and $L \in \mathbb{R}^{n_z}$ such that the following approximation holds
\begin{equation}\label{eq:basic_approx}
B(\xi(q)) \mu \approx L^TZ
\end{equation}
while guaranteeing the following inequalities
\begin{equation}\label{eq:inc_cstr}
[\frac{dB}{d\eta}(\eta)] \mu \geq \epsilon \quad \forall \eta \in [0,1]
\end{equation}
given a positive constant $\epsilon$.

The constraint \eqref{eq:inc_cstr} expresses the fact that $\Gamma(\cdot)$ has to be strictly increasing in order to guarantee the existence of the inverse map $\Gamma^{-1}(\cdot)$. This constraint can be transformed into a finite number of linear inequalities by defining a sufficiently dense grid of $\eta$ over the interval $[0,1]$ as
\begin{equation}
\boldsymbol{\eta}=\{0=\eta_1 < \eta_2 \dots < \eta_{n_\text{grid}}=1\}
\end{equation}
and can be rewritten in the following compact form 
\begin{equation}\label{eq:mono_cstr_grid}
\underbrace{\begin{bmatrix}
\frac{dB^{(1)}}{d\eta}(\eta_1) & \dots  & \frac{dB^{(n_b)}}{d\eta}(\eta_1) \\
\frac{dB^{(1)}}{d\eta}(\eta_2) & \dots  & \frac{dB^{(n_b)}}{d\eta}(\eta_2) \\
\vdots & \ddots & \vdots \\
\frac{dB^{(1)}}{d\eta}(\eta_{n_\text{grid}}) & \dots  & \frac{dB^{(n_b)}}{d\eta}(\eta_{n_\text{grid}})
\end{bmatrix}}_{[dB {}/ d\eta {} (\boldsymbol{\eta})]} \cdot
\underbrace{\begin{bmatrix}\mu_1\\\mu_2\\\vdots\\\mu_{n_b}\end{bmatrix}}_{\mu} 
\geq \epsilon \cdot \underbrace{\begin{bmatrix}1\\1\\\vdots\\1\end{bmatrix}}_{\boldsymbol{1}}     
\end{equation}

To summarize, the optimal parameters $(\mu,L)$ can be obtained by solving the following linear program (LP):
\begin{subequations}\label{eq:lp}
\begin{align}
& \underset{\mu,L}{\min} \quad \underset{(q,Z) \in \mathcal{D}}{\max} \omega(q,Z) \cdot |B(\xi(q)) \mu - Z^TL| \label{eq:cost_lp}\\
& \text{s.t. } \quad [\frac{dB}{d\eta}(\boldsymbol{\eta})] \mu \geq \epsilon \cdot \boldsymbol{1} 
\label{eq:mono_cstr}
\end{align}
\end{subequations}
where the \emph{weight indicator} $\omega(q,Z)$: $\mathbb{R} \times \mathbb{R}^{n_z} \to \mathbb{R}_{>0}$ is introduced in the cost function \eqref{eq:cost_lp} to enforce specific precision. This formulation obviously recalls the known \emph{weighted norm approximation} where $\infty$-norm is employed. Alternative formulation, based on the $L_2$-norm, can also be adopted leading to a quadratic programming (QP) problem.

\subsection{Preliminary Analysis}\label{subsec:pre_analysis}

For convenience, the quadruple \emph{identification parameters} $(n_m,\beta,\epsilon,\boldsymbol{\eta})$ is defined as follows:
\begin{itemize}
\item The parameter $n_m \geq 1$ is used to manipulate the degree of freedom;
\item The parameter $\beta$ is normally fixed as $0.5$;
\item The parameter $\epsilon$ is normally fixed as $1$;
\item The output grid $\boldsymbol{\eta}$ need to be sufficiently dense.
\end{itemize}

The strictly positive weight indicator $\omega(q,Z)$ can be  small or large over some subsets of data while being equal to $1$ for the remaining. The strict positivity of this weight is to avoid any drastically loose approximation as well as to guarantee the well-posedness of the formulation. In this paper, we utilize a simple form of the weight indicator $\omega(q,Z)$ as 
\begin{equation}
\omega(q,Z) =
\begin{cases}
\rho_1 & \text{if } (q,Z) \in \mathcal{W}_1 \\
\rho_2 & \text{if } (q,Z) \in \mathcal{W}_2 \\
\thinspace \vdots \\
1 & \text{otherwise}
\end{cases}
\end{equation}
where $\rho_i > 0$ are the constant weights corresponding to disjoint subspaces $\mathcal{W}_i \subset [\underline{q}, \overline{q}] \times \mathcal{Z}$.

Once the identification parameters are set and the weight indicator is constructed, the \emph{model parameters} $(\mu,L)$ can be obtained owing to the following feasibility assessment:

\begin{proposition}\label{prop:id_feas}
The LP \eqref{eq:lp} is feasible.
\end{proposition}
\textbf{Proof}. According to \eqref{eq:basis_fcn} and \eqref{eq:mono_cstr_grid}, $[dB {}/ d\eta(\boldsymbol{\eta})]$ is an $n_\text{grid} \times n_b$ matrix with the elements of the first column being zeros while the remaining being strictly positive. Hence, it is obvious that there exists $\mu$, for instance, the one with $\mu_i > 0$ for all $i=2,\dots,n_b$, such that the l.h.s. of \eqref{eq:mono_cstr} is strictly positive. Thus, the fulfillment of the feasibility condition \eqref{eq:mono_cstr} is guaranteed with any $\epsilon > 0$. $\hfill \Box$ 

The identification residual of the identified model can be characterized as follows:

\begin{proposition}\label{prop:id_residual}
Given the map $F$ whose model parameter, namely $(\mu,L)$, is a feasible solution of the LP \eqref{eq:lp} with corresponding minimization cost
\begin{equation}
J = \underset{(q,Z) \in \mathcal{D}}{\max} \omega(q,Z) \cdot \left|B(\xi(q)) \mu - Z^TL\right|
\end{equation}
The identification residual at any learning data point $(q,Z)$ is such that:
\begin{equation}\label{eq:id_residual}
\frac{|q-F(Z)|}{\overline{q}-\underline{q}} \leq \frac{1}{\omega(q,Z)} \cdot \frac{J}{\epsilon}
\end{equation}
\end{proposition}
\textbf{Proof}. Consider any $(q,Z) \in \mathcal{D}$. Let us denote the approximate as $\hat{q}=F(Z)$, i.e. $B(\xi(\hat{q}))\mu=L^TZ$. Since the map $\Gamma(\cdot)$ is continuous differentiable, one clearly has
\begin{equation*}
\begin{aligned}
J & \geq \omega(q,Z) \cdot |B(\xi(q))\mu - L^TZ|\\
& = \omega(q,Z) \cdot |B(\xi(q))\mu-B(\xi(\hat q))\mu|\\
& \geq \underset{\xi \in [0,1]}{\min} \left[[\frac{dB}{d\xi}(\xi)]\mu\right] \omega(q,Z) \cdot |\xi(\hat{q})-\xi(q)|\\
& \geq \epsilon \cdot \omega(q,Z) \cdot \frac{|\hat{q}-q|}{\overline{q}-\underline{q}} 
\end{aligned}
\end{equation*}
which is equivalent to \eqref{eq:id_residual}. $\hfill \Box$ 

Proposition \ref{prop:id_residual} implies that the desired fit would be obtained for a sufficiently small minimization cost. In other words, the approximations \eqref{eq:approx_map} and \eqref{eq:basic_approx} are relevant.

\section{Piecewise Nonlinear System Identification}\label{sec:pwnl}  

The lack of universal property of the approximator \eqref{eq:approx_struct} need to be underlined. In general, finding $L$ and $\Gamma(\cdot)$ such that the approximation $\Gamma(q) \approx L^TZ$ holds might be impossible for many identification problems even for infinite number of parameters. This structural limitation can be overcame by using a number of submodels where each submodel valids over a region of the regression domain. An identification methodology, which is a generalization of piecewise affine frameworks , is proposed in this section to do so.

\subsection{Problem Statement and Proposed Methodology}

\begin{problem}\label{pb:pwma}
Given the data $\mathcal{D}=\{(q(k),Z(k)\}_{k=1}^N$ where $Z \in \mathcal{Z} \subseteq \mathbb{R}^{n_z}$ is the regressor and $q \in [\underline{q},\overline{q}] \subset \mathbb{R}$ is the output. Find $s$ maps $F_{(i)} : \mathbb{R}^{n_z}  \rightarrow \mathbb{R}$ of the form: 
\begin{equation}
 	F_{(i)}(Z):=\Gamma_{(i)}^{-1}(L_{(i)}^TZ)
\end{equation}  
where $L_{(i)} \in \mathbb{R}^{n_z}$, the strictly increasing map $\Gamma_{(i)}(\cdot)$ is parametrized as
\begin{equation}
 	\Gamma_{(i)}(q)=B(\xi(q))\mu_{(i)} \text{; } \xi(q)=\frac{q-\underline{q}}{\overline{q}-\underline{q}}
\end{equation} 
and corresponding regions $\{ \mathcal{R}_{(i)} \}_{i = 1}^s$ forming a complete partition of $\mathcal{Z}$ (i.e. $\bigcup_{i = 1}^s \mathcal{R}_{(i)}=\mathcal{Z}$ and $\mathcal{R}_{(i)} \cap \mathcal{R}_{(j)}=\emptyset \thinspace \forall i \neq j$) such that the following approximation holds:
\begin{equation}\label{eq:pb_pwma}
 q \approx
  \begin{cases}
   F_{(1)}(Z) & \text{if } Z \in \mathcal{R}_{(1)} \\
   \qquad \vdots \\
   F_{(s)}(Z) & \text{if } Z \in \mathcal{R}_{(s)} \\
  \end{cases}
\end{equation}
where $s$ is referred to as the model complexity.
\end{problem}

Once the partition $\{ \mathcal{R}_{(i)} \}_{i = 1}^s$ is constructed, the corresponding classification $\{ \mathcal{D}_{(i)} \}_{i = 1}^s$ can be defined as
\begin{equation}\label{eq:dataset_region}
\mathcal{D}_{(i)}=\{ (q,Z) \in \mathcal{D}: Z \in \mathcal{R}_{(i)} \} 
\end{equation}
According to Proposition \ref{prop:id_feas}, the \emph{submodel parameters} $(\mu_{(i)},L_{(i)})$ of $F_{(i)}(\cdot)$ are hereafter available by solving the following LP:
\begin{subequations}\label{eq:lp_i}
\begin{align}
& \underset{\mu_{(i)},L_{(i)}}{\min} \quad \underset{(q,Z) \in \mathcal{D}_{(i)}}{\max} \omega(q,Z) |B(\xi(q)) \mu_{(i)} - Z^TL_{(i)}|\\
& \text{s.t. } \quad [\frac{dB}{d\eta}(\boldsymbol{\eta})] \mu_{(i)} \geq \epsilon \cdot \boldsymbol{1} 
\end{align}
\end{subequations}

The challenge is to find the minimal number $s$ of regions $\mathcal{R}_{(i)}$ such that the following conditions hold:
\begin{enumerate}
\item The associated dataset is such that 
\begin{equation}\label{eq:nonempty}
\#\mathcal{D}_{(i)} \geq \kappa \quad \forall i=1,\dots,s
\end{equation}
where $\kappa$ is the minimum required cardinality of each dataset. This lower bound is imposed to avoid overfit.
\item The associated identification residual is bounded by a given tolerance $\sigma > 0$, i.e.
\begin{equation}\label{eq:fit_bound}
\gamma_{(i)}=\underset{(q,Z) \in \mathcal{D}_{(i)}}{\max} \frac{|q-F_{(i)}(Z)|}{\overline{q}-\underline{q}} \leq \sigma \quad \forall i=1,\dots,s
\end{equation} 
\end{enumerate}

The structure \eqref{eq:pb_pwma} satisfying \eqref{eq:fit_bound} will be reffered to as an \emph{$s$-regions $\sigma$-error} PWNL model. The following heuristic identification methodology is proposed to derive such a model, consisting of two stages:

\begin{enumerate}
\item \emph{Domain partitioning:} At this stage, the partition is initially rough and is then refined at each step. Precisely, we seeks partition $\{ \mathcal{R}_{(i)} \}_{i = 1}^r$ and associated classification $\{ \mathcal{D}_{(i)} \}_{i = 1}^r$, leading to a $r$-regions $\sigma$-error approximation with a \emph{raw complexity} $r$ for the desired fit. Generally, the lower $\sigma$, the finer resulting partition.
\item \emph{Complexity reduction:} At this stage, the initial partition is the one resulted from the previous stage and is simplified at each step. Precisely, we aim to obtain an approximation with the \emph{reduced complexity} $s < r$ by merging iteratively regions and datasets while maintaining the fitting level of $\sigma$-error.  
\end{enumerate}

The following sub-sections detail such stages.

\subsection{Domain Partitioning}\label{subsec:domain_partitioning}

For convenience, all involved variables should be normalized. Hence, we restrict our attention to hyperrectangular regression domains $\mathcal{Z} = [0, 1] \times \dots \times [0, 1] \subset \mathbb{R}^{n_z}$ without loss of generality.

An iterative algorithm is proposed for the nonuniform hyperrectangular partitioning of PWNL function domains, i.e. each region $\mathcal{R}_{(i)}$ is a hyperrectangle, satisfying both \eqref{eq:nonempty} and \eqref{eq:fit_bound}. In addition, we utilize $\nu$ as the minimum allowed size of the hyperrectangles along each dimension. This lower bound is not only imposed by the machine precision but also can be tuned by the designer. The algorithm is summarized as follows:

\begin{enumerate}
\item \emph{Initialization:} The algorithm starts with a rough initial partition which can be the whole domain, i.e. $r^{[0]}=1$, $\mathcal{R}_{(1)}^{[0]}=\mathcal{Z}$ and $\mathcal{D}_{(1)}^{[0]}=\mathcal{D}$.  
\item \emph{Algorithm iterations:} At each iteration $k$, a hyperrectangle $\mathcal{R}_{(\cdot)}^{[k]}$ with $\gamma_{(\cdot)}^{[k]} > \sigma$ is split into two subhyperrectangles by an axis-orthogonal hyperplane that goes through the centroid of the dataset $\mathcal{D}_{(\cdot)}^{[k]}$. The basic idea is to select the hyperplane leading to a maximal reduction in the identification residual while guaranteeing the size of the resulted subhyperrectangles and the cardinality of the associated datasets are not lower than $\nu$ and $\kappa$, respectively. In the sequel, the selected hyperplane will be replaced by a maximum-margin hyperplane that is parallel with it. Then, a new iteration of the algorithm is performed by increasing the complexity by one unit $r^{[k+1]}=r^{[k]}+1$. Otherwise, the split fails when resulted subhyperrectangles and subdatasets are too small for all split options.
\item \emph{Termination:} The algorithm ends when the identification residual is not larger than $\sigma$ or in the case of splitting failure. 
\end{enumerate}

The algorithm successfully terminates when, at the end of the iterations, a $\sigma$-error approximation is achieved while the size of all hyperrectangles and the cardinality of all datasets are not lower than $\nu$ and $\kappa$, respectively.

\subsection{Complexity Reduction}\label{subsec:complex_reduc}

Let us assume that the domain partitioning stage returns a partition $\{ \mathcal{R}_{(i)} \}_{i = 1}^r$, a classification $\{ \mathcal{D}_{(i)} \}_{i = 1}^r$ and the raw complexity $r$. If the model complexity $r$ needed for data fit is overestimated, this number can be reduced by forcing the parameters of some partitions to be identical while maintaining the fitting level in a postprocessing step as follows:

\begin{enumerate}
\item \emph{Initialization:} The algorithm starts with an admissible initial partition $\{ \mathcal{R}_{(i)} \}_{i = 1}^r$, i.e. $s^{[0]}=r$.  
\item \emph{Algorithm iterations:} At each iteration $k$, a nonempty subset $\mathcal{I} \subseteq \{ 1,\dots,s^{[k]} \}$ is computed satisfying the following conditions:
\begin{subequations}
\begin{align}
& \begin{bmatrix}\mu_{(i)}\\L_{(i)}\end{bmatrix}=\begin{bmatrix}\mu_{(j)}\\L_{(j)}\end{bmatrix} & \forall i \neq j \text{; } i,j \in \mathcal{I} \label{eq:compatible_subset1}\\
& \gamma_{(i)} \leq \sigma & \forall i \in \mathcal{I} \label{eq:compatible_subset2}
\end{align}
\end{subequations}
Then, a new iteration is performed by reducing the model complexity $s^{[k+1]}=s^{[k]}-\#\mathcal{I}+1$. Moreover, one merges all the regions and datasets whose indexes belonging to $\mathcal{I}$.
\item \emph{Termination:} The algorithm always successfully terminates and stops when the maximum number of iterations $N_\text{iter}$ is reached.
\end{enumerate}

The search for the subset $\mathcal{I}$ at each iteration $k$ will be detailed. Since searching for the subset $\mathcal{I}$ over the whole set of regions $\{ 1,\dots,s^{[k]} \}$ is computationally expensive, the search over several regions with indexes belonging to a set $\mathcal{L} \subseteq \{ 1,\dots,s^{[k]} \}$ might be more practical. Once the set $\mathcal{L}$ is chosen, its subset $\mathcal{I}$ satisfying \eqref{eq:compatible_subset1} and \eqref{eq:compatible_subset2} is referred to as a \emph{compatible subset}. The searching for the compatible subset $\mathcal{I}$ is summarized as follows:

\begin{enumerate}
\item \emph{Initialization:} The procedure start with an initial guess $\mathcal{I}^{[0]}=\mathcal{L}$.
\item \emph{Algorithm iterations:} At each sub-iteration $l$, consider the following LP whose decision variables are $\{ \mu_{(i)},L_{(i)} \}_{i \in \mathcal{L}}$ and $\tilde{\mu},\tilde{L}$: 
\begin{subequations}\label{eq:mpLP}
\begin{align}
& \min \underset{i \in \mathcal{L}}{\sum} \thinspace \underset{(q,Z) \in \mathcal{D}_{(i)}}{\max} \omega(q,Z) |B(\xi(q)) \mu_{(i)} - Z^TL_{(i)}| \label{eq:mpLP_cost}\\
& \text{s.t. } \quad [\frac{dB}{d\eta}(\boldsymbol{\eta})] \mu_{(i)} \geq \epsilon \cdot \boldsymbol{1} \quad \forall i \in \mathcal{L} \label{eq:mpLP_mono_cstr}\\
& \text{and } \quad ||\begin{bmatrix}\mu_{(i)}\\L_{(i)}\end{bmatrix}-\begin{bmatrix}\tilde{\mu}\\\tilde{L}\end{bmatrix}||_\infty=0 \quad \forall i \in \mathcal{I}^{[l]} \label{eq:mpLP_comp_cstr}
\end{align}
\end{subequations}
The feasibility of the LP \eqref{eq:mpLP} is straightforward from Proposition \ref{prop:id_feas}. The additional constraints \eqref{eq:mpLP_comp_cstr} are referred to as \emph{compatibility constraints}, including $\#\mathcal{I}^{[l]}$ linear equalities. If the LP \eqref{eq:mpLP} results into the submodels with the following identification residual
\begin{equation}\label{eq:residual_merge}
\gamma_{(\mathcal{I}^{[l]})}=\underset{i \in \mathcal{I}^{[l]}}{\max} \thinspace \gamma_{(i)}
\end{equation}
which is larger than $\sigma$, it is an indication that the minimization cost of \eqref{eq:mpLP} is insufficiently small. This cost can be reduced by removing some of the compatibility constraints \eqref{eq:mpLP_comp_cstr}. Since the sensitivities of \eqref{eq:mpLP} with respect to such constraints are represented by the corresponding Lagrange multipliers, namely $\{ \lambda_{(i)} \}_{i \in \mathcal{I}}$, a common strategy is to iteratively remove the constraint associated with the largest Lagrange multiplier. This is equivalent to removing an element of $\mathcal{I}$, i.e.
\begin{equation}
\mathcal{I}^{[l+1]} = \mathcal{I}^{[l]} \setminus \{\underset{i \in \mathcal{I}^{[l]}}{\text{argmax}} \thinspace \lambda_{(i)}\}
\end{equation}
\item \emph{Termination:} The procedure stops when $\gamma_{(\mathcal{I}^{[l]})} \leq \sigma$.
\end{enumerate}

\section{Approximate Explicit Model Predictive Control}\label{sec:mpc} 

\subsection{Extension to Multi-Input Multi-Output MPC}

Although the presented methodology is restricted to the multi-input single-output relationships, the extension to approximate an MPC control law is trivial by approximating each component of the control vector (see Sub-section \ref{subsec:wcs} for example).

\subsection{On the Stabilization}

A crucial aspect of the resulting EMPC in the regulation problem is setpoints stabilization. On the one hand, the key properties to guarantee stabilization of such an approximate control law have been extensively studied in the literature. Example is \cite{grancharova2012explicit} which indicates the boundedness of the state trajectory and its convergence to a neighborhood of the setpoint can be achieved by a sufficient small approximation error. A quantitative analysis on such an upper bound is presented in \cite{canale2010efficient}, however, is difficult to applied here due to the limitations of the identification approach. On the other hand, the efforts to establish asymptotic stability in a setpoint's neighborhood generally deploy the dual-mode design where a local stabilizing controller is used. In this paper, we impose a sufficiently high weight for the data point in a neighborhood of the setpoint to acquire an extremely high local accuracy. This parameter choice is based on the approximate stability result of discontinuous control laws (Theorem 1, Assumption 3, \cite{scokaert1999suboptimal}).

\subsection{Remarks on Weight Tuning}

This can be generalized to the reference tracking problem. First, we extract from the learning data a discrete set of samples corresponding to stationary state of the closed loop. Second, in the identification procedure of control actuators, a typically high weight, namely $\rho_\text{st}$, is associated to such stationary data to obtain nearly zero steady errors of the closed-loop controlled by the EMPCs. Lastly, a moderate weight, namely $\rho_\text{stab}$, is associated to the data points in the neighborhood of the considered stationary data.  

Note finally that a piecewise control law may results into undesired oscillations in the closed-loop response. This can be partially overcame by imposing a moderate weight, namely $\rho_\text{sw}$, for the data in the neighborhoods of the switching surfaces between computed regions.

\section{Illustrative Examples}\label{sec:apps}

\subsection{Application 1: Reference Tracking Control of A Stirling Engine Based Power Generation System}\label{subsec:stirling}

\begin{figure}
\centering
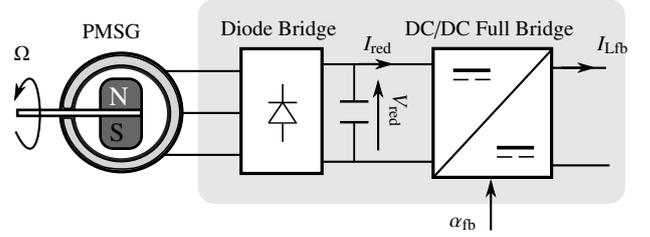
\caption{Schematic diagram of a power electronics module in thermodynamic power plants. The component named PMSG is permanent magnetic synchronous generator whose shaft is connected with a Stirling engine.}
\label{fig:stirling_diagram}
\end{figure}

Consider the following system which belongs to a stand-alone hybrid Stirling engine / supercapacitor power plant (see Fig.~\ref{fig:stirling_diagram}) with $x_1=\Omega$ being the shaft rotational speed, $x_2=I_\text{red}$ and $x_3=V_\text{red}$ being the rectified current and the rectified voltage while $x_4=I_\text{Lfb}$ and $u=\alpha_\text{fb}$ being the output current and the duty ratio of the DC/DC full bridge converter:
\begin{equation}\label{eq:stirling_sys}
\begin{array}{lcl}
\dot{x}_1 & = & -a_1x_1-a_3x_2+a_2\\
\dot{x}_2 & = & -a_4x_2+a_6x_1-a_7x_3\\
\dot{x}_3 & = & a_8 (x_2-kx_4u)\\
\dot{x}_4 & = & a_9 (-x_5^\text{st}+kx_3u)
\end{array}
\end{equation}
The parameters of this model are $a_1 = − 0.183$, $a_2 = 558.11$, $a_3 = 118.4453$, $a_4 = 9615.4$, $a_6 = 5101.1$, $a_7 = 641.02$, $a_8 = 425.53$, $a_9 = 6666.7$, $k = 0.5$ and $x_5^\text{st} = 50$. The system states should satisfy positivity constraints $x \geq 0$ since the used converter is not reversible. The control variable $u$ is the duty ratio of this converter and has a strong saturation $u\in [0,~1]$. The sampling time is fixed as $T_s=100\mu$s.

The control objective is to force the state $x_4$ to track a reference signal $x^r_4$ while respecting the constraints. In \cite{rahmani2013control}, it has been shown that this is a challenging problem as \eqref{eq:stirling_sys} contains highly oscillatory modes that induce constraints violation if the latter are not explicitly addressed. These make the NMPC methodology particularly suitable. This can be achieved by minimizing the following cost:
\begin{equation}
\begin{aligned}
\sum_{i=1}^{N_p} ( || x(k+i) - x^r(k) ||^2_Q + || u(k+i-1) - u^r(k) ||^2_R ) + \rho \delta^2
\end{aligned}
\end{equation}
with the prediction horizon $N_p = 3$, the cost matrices $Q = I$, $R = 1$ and the weight penalizing constraint violations $\rho = 10^4$ while fulfilling the hard constraint $0 \leq u(k) \leq 1$, the soft constraints $[0-\delta,~4.5-0.1 \delta,~55-\delta,~1-\delta]^T \leq x(k) \leq [40+\delta,~5.5+0.1 \delta,~200+\delta,~25+\delta]^T$ where $\delta \geq 0$. The prediction simulator is under-sampled with a sampling period of $\tau_u=1$ms.

The fast dynamics of the considered system would require the design of an EMPC controller. In order to do this, a data record is generated by simulating the closed-loop system with the implicit NMPC during a $10$s scenario. This simulation scenario is based on the considered range of the reference as $x^r_4 \in [7, 13]$ [A]. Hereafter, a learning data is extracted with the cardinality of $N = 24205$ including $48$ stationary data points that are then used to build the instances
\begin{equation}
Z(k)=\begin{bmatrix}x^T(k), \left(x^r(k)\right)^T, \left(u^r(k)\right)^T, \left(e(k)\right)^T\end{bmatrix}^T; q(k) = u(k)
\end{equation}
The reference tracking error $e(k)=x(k)-x^r(k)$ is incorporated into the regressor $Z(k)$ to further reduce the partition complexity of the explicit controller. The values $\underline{q}=0.35$, $\overline{q}=0.7$ are the range of $u$ in the learning data. 

Several approximate controllers have been identified  in Matlab using the IBM's CPLEX with the same identification and partitioning parameters. The identification parameters have been chosen with $\beta=0.5$, $\epsilon=1$, $n_m=10$ and a uniform grid $\boldsymbol{\eta}$ of $n_\text{grid}=50$ points. The weight indicator is defined by imposing a high weight of $\rho_{\text{st}}=10^2$ for stationary data, a moderate weight of $\rho_{\text{stab}}=10$ for data in the neighborhood of stationary data, a lower weight of $\rho_{\text{sw}}=2$ for data near the boundary of each hypercubes and the weight of $1$ for the remaining. The partitioning parameters have been chosen as $\nu=5 \times 10^{-2}$ and $\kappa=50$.

\begin{figure}
\centering
\includegraphics{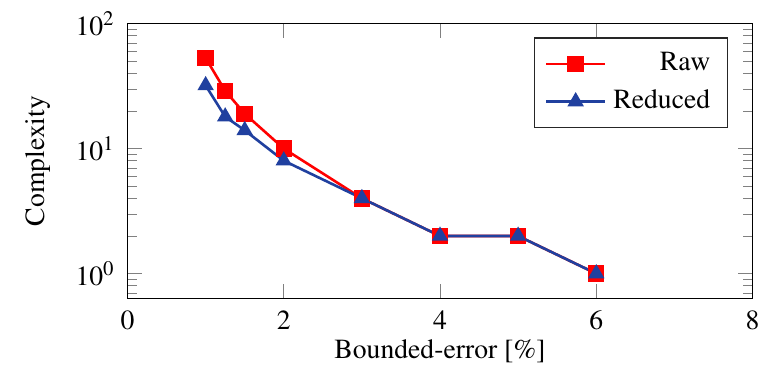}
\caption{Trade-off curve.}
\label{fig:stirling_empc_complex}
\end{figure}

\begin{figure}
\includegraphics{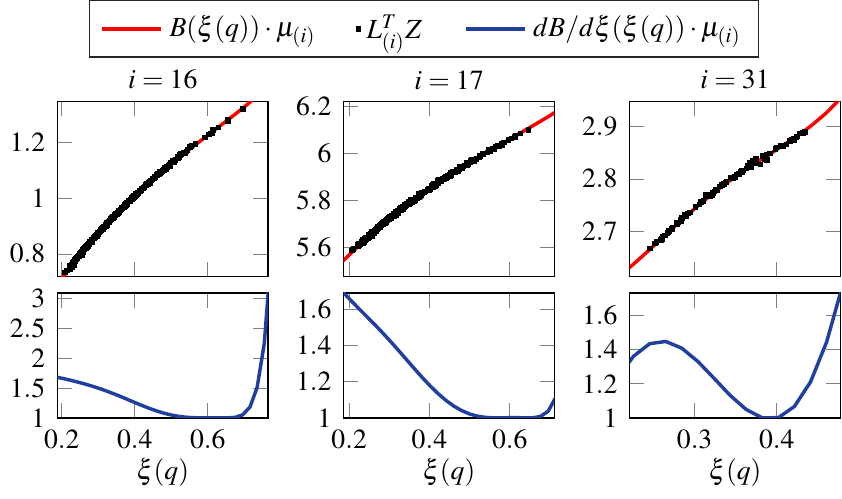}
\caption{Nonlinearity and strict monotonicity of the $i$-th submodels corresponding to the $32$-regions $1\%$-error EMPC. The black dots are the function evaluation for validation data while the red and the blue curves are the identified nonlinear function and its derivative, respectively.}
\label{fig:stirling_empc_static}
\end{figure}

Fig.~\ref{fig:stirling_empc_complex} portrays a complexity assessment of identified controllers with different values of the bounded-error $\sigma$ while the nonlinear characteristic of several $i$-th submodels corresponding to a $32$-regions $1\%$-error EMPC are illustrated in Fig.~\ref{fig:stirling_empc_static}. In Fig.~\ref{fig:stirling_empc_static}, the red curves represent the identified maps $\Gamma_{(i)}^{-1}(\xi)$ which are strictly increasing as their derivatives (blue curves) are higher than $\epsilon=1$. Although these maps are defined over the interval $\xi \in [ 0,1 ]$, we only show the range with the presence of data points. The black dots illustrates the nonlinear nature of data.

\begin{figure}
\centering
\includegraphics{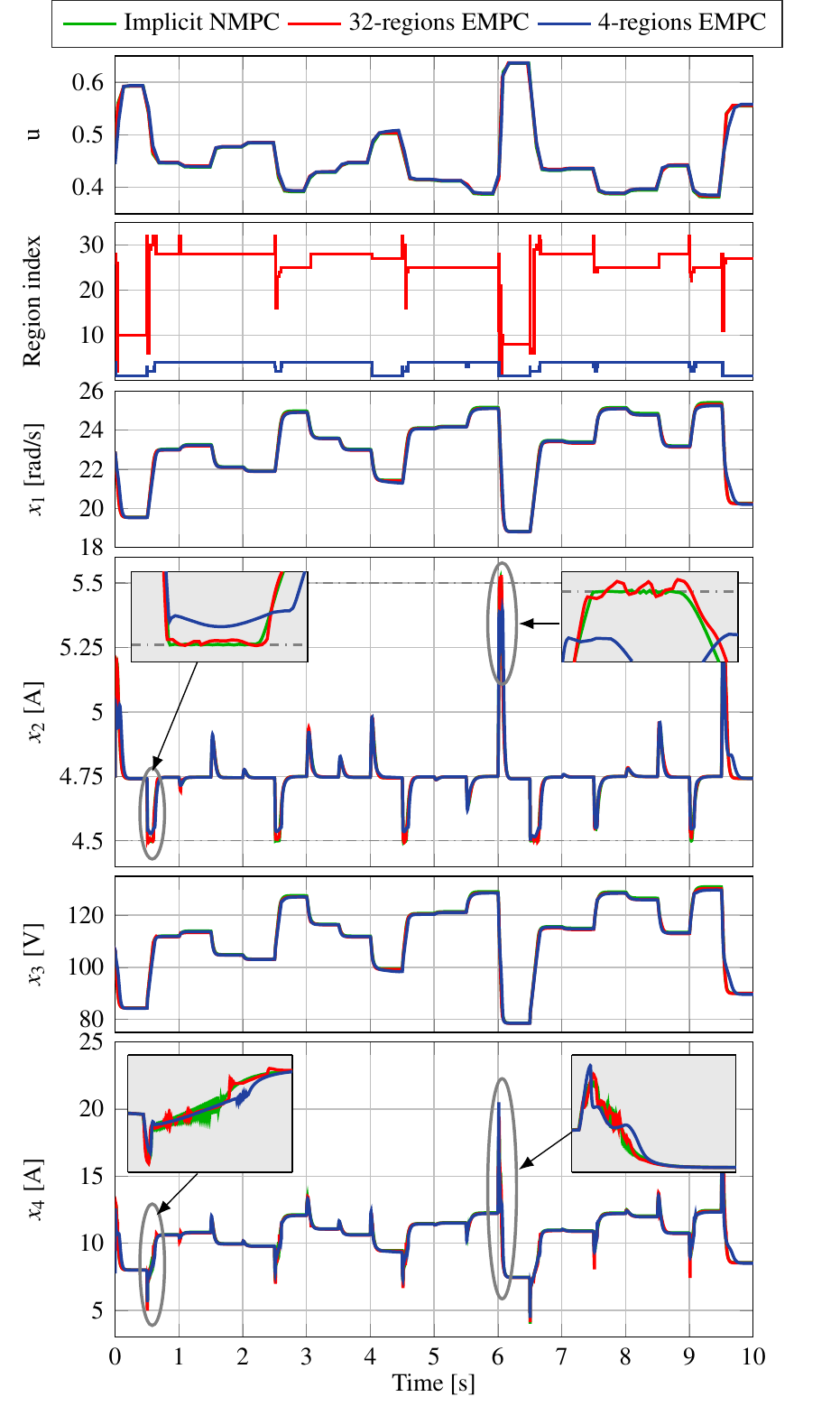} 
\caption{Closed-loop performance comparison. The constraint satisfaction are portrayed in the zooms on the evolution of $x_2$.}
\label{fig:stirling_empc_valid}
\end{figure}

The performance of the resulted EMPCs is evaluated by simulations and is compared with the closed-loop controlled by the implicit NMPC. Fig.~\ref{fig:stirling_empc_valid} shows the closed-loop response with different controllers: the implicit NMPC (green), the $32$-regions $1\%$-error EMPC (red) and the $4$-regions $3\%$-error EMPC (blue). Notably, there are switching between regions, especially in the transient duration. It is also clear that the tracking errors of those EMPCs are negligible.

\begin{figure}
\centering
\includegraphics{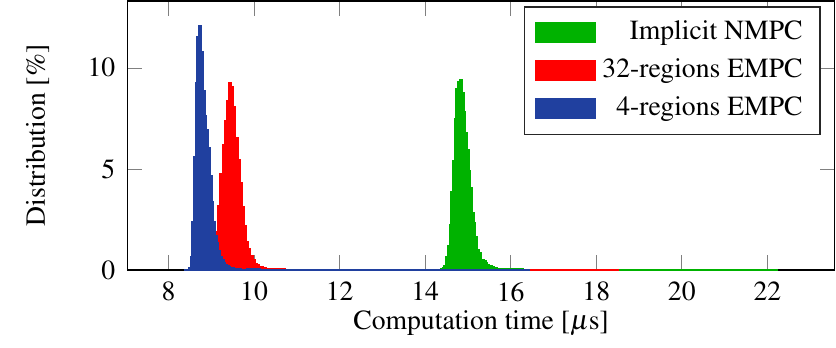}
\caption{Computation time. The implicit NMPC was built using the ACADO code generation \cite{quirynen2015autogenerating} with the qpOASES solver \cite{ferreau2014qpoases} while the EMPCs were implemented in C language. Platform: $2.6$ GHz Intel(R) Core(TM) i$7$ and $16$GB of RAM.}
\label{fig:stirling_empc_time}
\end{figure}

Regarding real-time aspect, Fig.~\ref{fig:stirling_empc_time} illustrates the computation efficiency of the EMPCs. Their computation time vary between $8$$\mu$s and $18$$\mu$s which is slightly faster than that of the implicit NMPC delivered by ACADO and obviously suitable for this application.

\subsection{Application 2: Control of A Compression Station}\label{subsec:wcs}

\begin{figure}
\begin{center}
\includegraphics[scale=0.9]{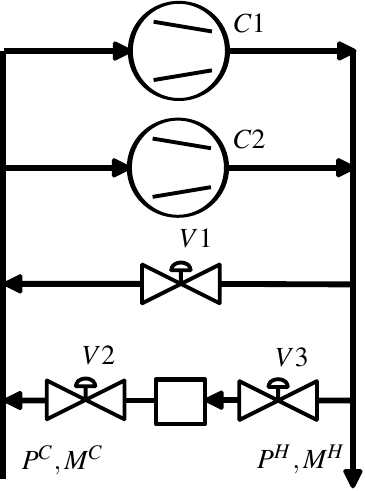}
\caption{Schematic diagram of a warm compression station. Arrows represents the fluid's direction. The block between $V2$ and $V3$ valves represents the capacity buffer.}
\label{fig:wcs_sbt}
\end{center}
\end{figure}

To illustrate the efficiency of the proposed framework for a multi-input multi-output case, we consider the constrained control for a warm compression station \cite{bonne2016wcs} (see Fig. \ref{fig:wcs_sbt}). Variables to be regulated are the low pressure $P^C$ and the high pressure $P^H$. The unmeasured flows $M^H$ and $M^C$, respectively denote the inflow and outflow, are considered as the process disturbances. $M_{C1}$ and $M_{C2}$ respectively denote the flowrates handled by the first and second compressor which are variable frequency driven, $M_{V1}$ represents the flowrate passing through the bypass valve V1 while $M_{V2}$ and $M_{V3}$ depict the flows generated by the charge valve V2 and discharge valve V3. These five variables are the manipulated variables of the process. Valves opening is obviously in the range $[0, 100]$ [\%], the compressor C1 is limited to a speed in the range $[30, 53]$ [Hz] while the compressor C2 is manually set at a speed of $50$ [Hz]. A simulator built up from the SimCryogenics Toolbox \cite{bonne2015simulink} is available.

In \cite{bonne2016wcs}, the three variables $M_{prod}=M_{C1}+M_{C2}+M_{V1}$, $M_{ch}=M_{V2}$ and $M_{dis}=M_{V3}$ are introduced as the virtual actuators. Then, the virtual actuators value (in [kg/s]) will have to be converted into actual control actions for the system actuators (in [Hz] or [\%]) by a conversion function. Then, a discrete-time linear model in the following form
\begin{equation}
    x(k+1)=A \cdot x(k)+B \cdot u(k)+F \cdot d(k)
\end{equation}
with the sampling time $T_s=0.25$s is synthesized for the control design. The system states, the control inputs and the disturbances are defined as the deviations around an operating point of interest, i.e. 
\begin{equation}
    x=\begin{bmatrix} \Delta P^H \\ \Delta P^C \end{bmatrix} \text{; } u=\begin{bmatrix} \Delta M_{prod} \\ \Delta M_{ch} \\ \Delta M_{dis} \end{bmatrix} \text{; } d=\begin{bmatrix} \Delta M^H \\ \Delta M^C \end{bmatrix}
\end{equation}

The control objective is to regulate system states to the setpoint $P^H=16$ [bar], $P^C=1.05$ [bar] under variable disturbances as well as input saturation. In addition, the low pressure is limited by a soft constraint $P^C \in [1, 1.1]$ [bar]. In \cite{bonne2016wcs}, a MPC design is derived with the following cost
\begin{equation}
\begin{aligned}
\sum_{i=1}^{N_p} ( || x(k+i) ||^2_Q + || u(k+i-1) - u^r(k) ||^2_R )
\end{aligned}
\end{equation}
in which the control reference $u^r$ is computed based on the estimated disturbances. A typically long prediction horizon of $N_p=200$ has been used.

Our objective is to identify three maps $F_i$ such that $q_i \approx F_i(Z)$, $i=1, 2, 3$ where: 
\begin{equation}
    Z(k)=\begin{bmatrix}x(k) \\ u^r(k)\end{bmatrix} \text{; } q_i(k) = u_i(k)
\end{equation}
A data record over an interval of $25 \times 10^3$s is available. It is hereafter extracted into a learning data with the cardinality of $N = 11633$ including $388$ stationary data points. The identification parameters have been chosen with $\beta=0.5$, $\epsilon=1$, $n_m=10$ and a nonuniform grid $\boldsymbol{\eta}$ of $n_\text{grid}=100$ points. The weight indicator is defined by $\rho_{\text{st}}=10^3$, $\rho_{\text{stab}}=20$ and $\rho_{\text{sw}}=5$. The partitioning parameters have been chosen as $\nu=1 \times 10^{-2}$ and $\kappa=10$.

\begin{figure}
\centering
\includegraphics{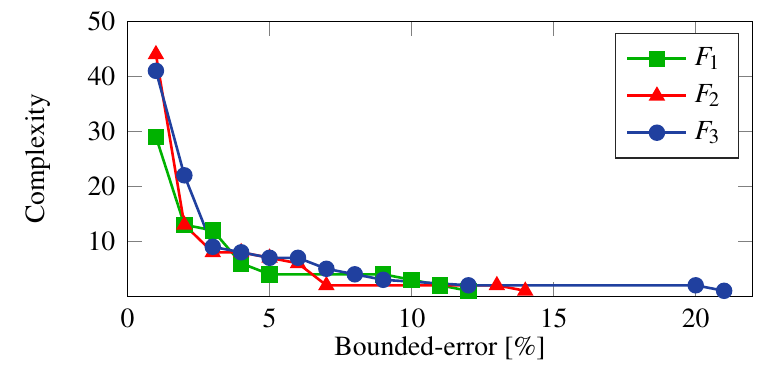}
\caption{Trade-off curve.}
\label{fig:wcs_empc_complex}
\end{figure}

Fig. \ref{fig:wcs_empc_complex} shows the trade off between the (reduced) complexity of the identified maps given different level of accuracy.

\begin{figure}[H]
\begin{subfigure}[b]{0.5\textwidth}
\centering
\includegraphics{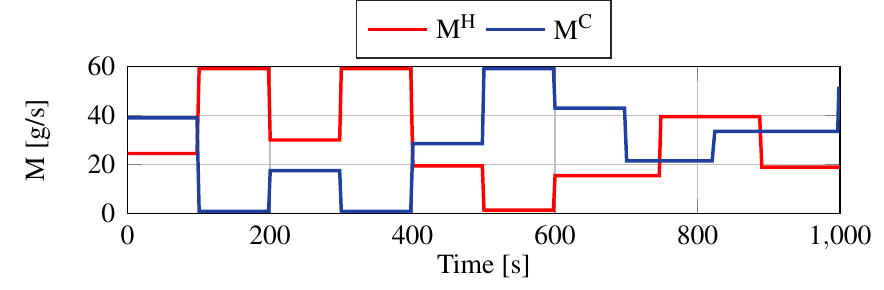}
\caption{Inflow and outflow}
\end{subfigure}
\begin{subfigure}[b]{0.5\textwidth}
\centering
\includegraphics{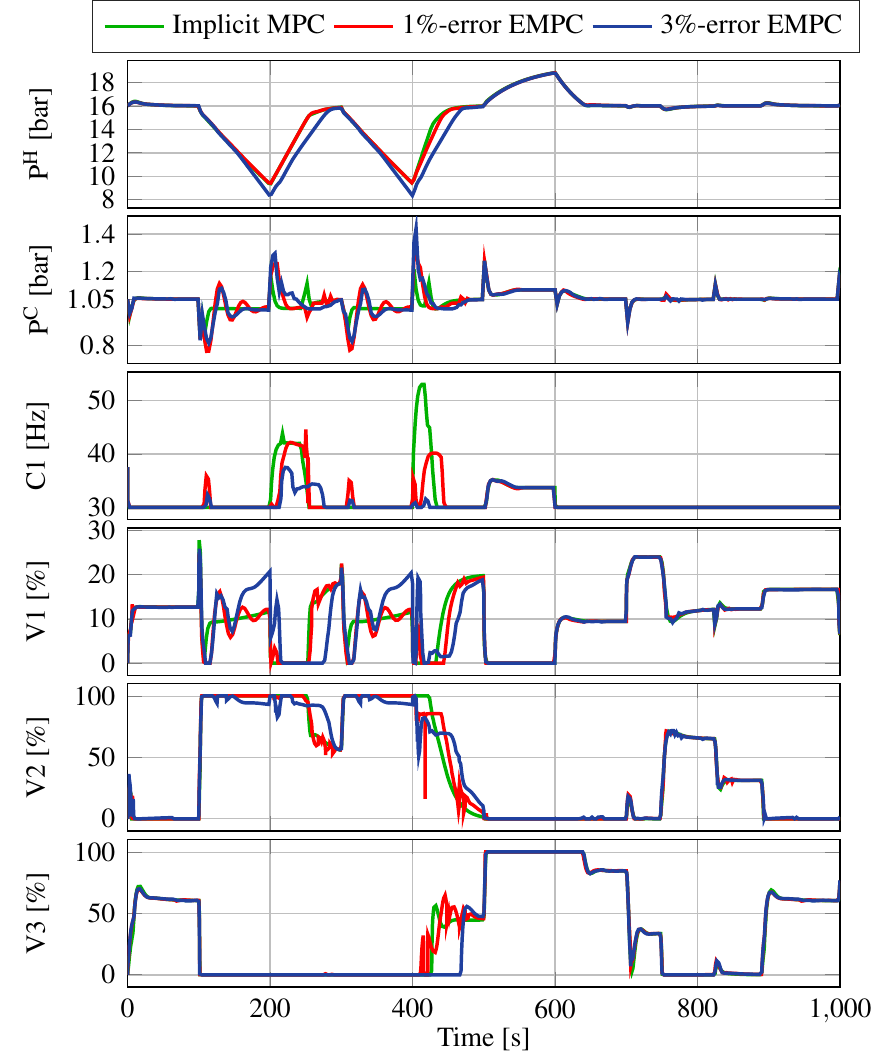}
\caption{Controlled pressures and manipulated actuators}
\end{subfigure}
\caption{Closed-loop performance comparison. The depicted response of the pressures has an amount of performance degradation.}
\label{fig:wcs_empc_valid}
\end{figure}

Fig. \ref{fig:wcs_empc_valid} depicts the closed-loop behaviours of the system in time under different controllers. The first $600$s duration shows how the system reacts when disturbed by some excessive difference between the inflow and the outflow. It is worth noting that pressure stability is ensure and large disturbances are effectively handled even in this critical situation. With less difference between inflow and outflow in the last $400$s duration, the performance of the implicit MPC and the identified EMPCs are obviously comparable. Notably, the pressures are returning to the setpoints without oscillations during such a normal situation.

\section{Conclusion and Future Work}\label{sec:conclusion}

In this paper, a new approach is proposed for explicit piecewise nonlinear presentation of approximate MPC control laws. Future research efforts will be devoted to different MPC settings such as output-feedback MPC, hybrid MPC or distributed MPC. Further developments include applying the methodology outside the MPC context, e.g. off-line computation of the moving
horizon observers.

\bibliographystyle{ieeetr}

\end{document}